\documentclass[aos,preprint,tikz,border=2pt,png]{imsart}

\setattribute{journal}{name}{}
\usepackage[hyphens]{url}
\RequirePackage[OT1]{fontenc}
\RequirePackage{amsthm,amsmath, amsfonts, amssymb, graphics, graphicx, subfigure, color, multirow}
\RequirePackage{natbib}
\RequirePackage[colorlinks,citecolor=blue,urlcolor=blue,breaklinks]{hyperref}
\usepackage{appendix}
\usepackage{mathtools}

\usepackage{algpseudocode}
\usepackage{algorithm}
\usepackage{comment}

\newtheorem{lem}{Lemma}
\newtheorem{cor}{Corollary}

\newtheorem{thm}{Theorem}
\newtheorem{prop}{Proposition}
\newtheorem{example}{Example}

\newtheorem{rmk}{Remark}
\newtheorem{definition}{Definition}

\newcommand{\iid}{\stackrel{\mathrm{iid}}{\sim}}

\newcommand{\bbR}{\mathbb{R}}

\newcommand{\vM}{\boldsymbol{M}}

\newcommand{\vm}{\boldsymbol{m}}

\newcommand{\bbE}{\mathbb{E}}
\newcommand{\bbP}{\mathbb{P}}

\newcommand{\vB}{\boldsymbol{B}}

\newcommand{\vD}{\boldsymbol{D}}
\newcommand{\vE}{\boldsymbol{E}}

\newcommand{\vI}{\mathbf{I}}

\newcommand{\vV}{\boldsymbol{V}}

\newcommand{\vx}{\boldsymbol{x}}
\newcommand{\vX}{\boldsymbol{X}}

\newcommand{\vy}{\boldsymbol{y}}

\newcommand{\vz}{\boldsymbol{z}}

\newcommand{\bbeta}{\boldsymbol{\beta}}

\newcommand{\bLambda}{\boldsymbol{\Lambda}}

\newcommand{\bSigma}{\boldsymbol{\Sigma}}

\newcommand{\bepsilon}{\boldsymbol{\epsilon}}

\newcommand{\boldeta}{\boldsymbol{\eta}}

\newcommand\independent{\protect\mathpalette{\protect\independenT}{\perp}}
\def\independenT#1#2{\mathrel{\rlap{$#1#2$}\mkern2mu{#1#2}}}

\newcommand{\epf}{\hfill $\Box$}

\begin{document}

\begin{frontmatter}

\title{Global Testing under Sparse Alternatives for Single Index Models}
\runtitle{SIM-Detection}
\thankstext{t1}{Zhao's research was supported in part by the NSF Grant IIS-1633283. Liu's research was supported in part by the NSF Grants DMS-1120368 and DMS-1613035.}
\thankstext{t2}{These authors contributed equally.}

\begin{aug}
\author{\fnms{Qian} \snm{Lin}\thanksref{m1}\ead[label=e1]{qianlin@mail.tsinghua.edu.cn}\thanksref{t1,t2}}
\author{\fnms{Zhigen} \snm{Zhao}\thanksref{m2}\ead[label=e2]{zhaozhg@temple.edu}\thanksref{t1,t2}}
\and    
\author{\fnms{Jun S.} \snm{Liu}\thanksref{m4}\ead[label=e4]{jliu@stat.harvard.edu}\thanksref{t1}}
 
   \runauthor{Q. Lin Z. Zhao and J. S. Liu}
\affiliation{
Tsinghua University\thanksmark{m1} 
Temple University\thanksmark{m2}  
Harvard University\thanksmark{m4} 
}
  
  \address{Qian Lin\\
    Center for Statistical Sciences\\
    and Department of Industrial Engineering\\
    Tsinghua University\\
    Beijing 100084 \\
    China\\
    \printead{e1}}    
 \address{Zhigen Zhao\\
    Department of Statistics \\
    Temple University\\
    342 Speakman Hall \\
    1810 N. 13th Street \\
    Philadelphia, PENNSYLVANIA, 19122 \\
    USA\\
    \printead{e2}}   
  \address{Jun S.  Liu\\
    Department of Statistics\\
    Harvard University\\
    1 Oxford Street \\
    Cambridge, MA 02138  \\
    USA\\
    Center for Statistical Science\\
    Tsinghua University\\
    Beijing\\
    China\\
    \printead{e4}}

\end{aug}

\begin{abstract}
For the single index model $y=f(\bbeta^{\tau}\vx,\epsilon)$ with Gaussian design, 
where $f$ is unknown and $\bbeta$ is a sparse $p$-dimensional unit vector with at most $s$ nonzero entries,  we are interested in testing the null hypothesis that $\bbeta$, when viewed as a whole vector, is zero against the alternative that some entries of $\bbeta$ is nonzero. Assuming  that $var(\bbE[\vx \mid y])$ is non-vanishing, 
we define the generalized signal-to-noise ratio (gSNR)  $\lambda$ of the model as  the unique non-zero eigenvalue of $var(\bbE[\vx \mid y])$. We show that if $s^{2}\log^2(p)\wedge p$ is of a smaller order of $n$, denoted as $s^{2}\log^2(p)\wedge p\prec n$, where $n$ is the sample size, one can detect the existence of signals if and only if  gSNR$\succ\frac{p^{1/2}}{n}\wedge \frac{s\log(p)}{n}$. 
Furthermore, if the noise is additive (i.e., $y=f(\bbeta^{\tau}\vx)+\epsilon$),
one can detect the existence of the signal if and only if  gSNR$\succ\frac{p^{1/2}}{n}\wedge \frac{s\log(p)}{n} \wedge \frac{1}{\sqrt{n}}$. It is rather surprising that the detection boundary for the single index model with additive noise 
matches that for linear regression models.
 These results pave the road for thorough theoretical analysis of single/multiple index models in high dimensions. 
\end{abstract}

\end{frontmatter}

\section{Introduction}

Testing whether a subset of covariates have any  relationship with a quantitative response is one of the central problems in statistical analyses. Most existing literature focuses on linear relationships. The analysis of variance (ANOVA) introduced by Fisher in 1920s has been a main tool for statistical analyses of experiments and routinely used in countless applications. Under  the simple normal mean model, 
\begin{align*}
y_{i}=\mu+\beta_{i}+\epsilon_{i}, \quad i=1,...,p,
\end{align*} 
where $\epsilon_{i}\stackrel{iid}{\sim }N(0,\sigma^2)$ and $\sum_i \beta_i=0$, one-way ANOVA tests the null hypothesis $H_0$ :  $\beta_{1}=...=\beta_{p}=0$ against the alternative hypothesis
$H_1$:  at least one $\beta_{j}\neq 0$.
Although the test cannot indicate which $\beta_j$'s are nonzero, ANOVA is powerful in testing the global null  against alternatives of the form $\{\bbeta \mid \sum_{j}\beta_{j}^{2}>r^{2}\}$. 
\cite{arias2011global} extended these results to the linear model 
\begin{align}
y=X\beta+\bepsilon,
\end{align}
where $\bepsilon\sim N(0, \sigma^2 \vI)$, to test whether all the $\beta_i$'s are zero.
This can be formulated as the following null and alternative hypotheses:
\begin{equation}\label{eqn:lr:detect}
\begin{aligned}
\begin{cases}
H_{0}:&\quad \beta_{1}=...=\beta_{p}=0\\
H_{s,r}:& \quad \bbeta\in \Theta_{s}(r)\triangleq \{\bbeta\in \bbR^{p}_{s}\mid \|\bbeta\|_{2}^{2}\geq r^{2}\}
\end{cases}
\end{aligned} 
\end{equation}
where 
$\bbR^{p}_{s}$ denotes the set of $s$-sparse vector in $\bbR^{p}$ with the number of non-zero entries being no greater than $s$. 
\cite{arias2011global} and \cite{ingster2010detection}  showed that  one can detect the signal if and only if  $r^{2}\succ\frac{s\log(p)}{n}\wedge \frac{p^{1/2}}{n}\wedge\frac{1}{\sqrt{n}}$. The upper bound is guaranteed by an asymptotically most powerful test based on higher criticism  \citep{donoho2004higher}.

The linearity or other functional form assumption  is often too restrictive  in practice. Theoretical and methodological developments beyond parametric models are important, urgent, and extremely challenging. As a first step towards nonparametric global testing, we here study the  single index model  $y=f(\bbeta^{\tau}\vx,\epsilon)$, where $f(\cdot)$ is an unknown function. Our goal is to test the global null hypothesis that all the $\beta_i$'s are zero. The first challenge is to find an appropriate  formulation of  alternative hypotheses because  $\|\bbeta\|_{2}^{2}$ used in (\ref{eqn:lr:detect}) is not even identifiable in single index models. 

When $rank(var(\bbE[\vx\mid y]))$ is nonzero in a single index model, the unique non-zero eigenvalue $\lambda$ of $var(\bbE[\vx \mid y])$ can be viewed as the generalized signal to noise ratio (gSNR) \citep{lin2018sparse}. In Section \ref{sec:description}, we show that for the linear regression model, this $\lambda$ is almost proportional to  $\|\bbeta\|_2$ when it is small. The alternative hypotheses in (\ref{eqn:lr:detect}) can be rewritten as $gSNR>r^2$. Because of this connection, we can treat $\lambda$ as the separation quantity for the single index model and consider the following contrasting hypotheses:
\begin{align*}
\begin{cases}
H_{0}:& \beta_{1}=\beta_{2}=...=\beta_{p}=0, \quad or \quad \mbox{gSNR}=0,\\
H_{a}:& \mbox{gSNR} \geq \lambda_0.
\end{cases}
\end{align*}
We show that, under certain regularity conditions, one can detect a non-zero gSNR if and only if $\lambda_0\succ \frac{s\log(p)}{n}\wedge \frac{p^{1/2}}{n}\wedge\frac{1}{\sqrt{n}}$ for the single index model with additive noise.

This is a strong and surprising result because this detection boundary  is the same as that for the linear model. Using the idea from the sliced inverse regression (SIR)  \citep{li1991sliced}, we  show that this boundary can be  achieved by the proposed {\bf S}pectral test {\bf S}tatistics using {\bf S}IR (SSS)  and SSS with ANOVA test assisted (SSSa).
Although SIR has been advocated as an effective alternative to linear multivariate analysis \citep{chen1998can}, the existing literature has not provided satisfactory theoretical foundations until recently \citep{lin2018optimal,lin2018sparse,lin2017consistency}. 
We believe that the results in this paper provide further  supporting evidence to the speculation  that ``SIR can be used to take the same role as linear regression in model building, residual analysis, regression diagnoses, etc" \citep{chen1998can}.

In Section \ref{sec:description}, after briefly reviewing the SIR and related results in linear regression, we state the optimal detection problem and a lower bound for single index models. In Section \ref{sec:test}, we first show that the correlation-based Higher Criticism (Cor-HC) developed for linear models fails  for single index models, and then propose a test to achieve the lower bound stated in Section \ref{sec:description}. Some numerical studies are included in Section \ref{sec:studies}. We list several interesting implications and future directions in Section \ref{sec:discuss}. Additional proofs and lemmas are included in Appendices.

\section{Generalized SNR for Single Index Models}\label{sec:description}

\subsection{Notation}
The following notations are adopted throughout the paper.
For a matrix $\vV$, we call the space generated by its column vectors the column space and denote it by $col(\vV)$. The $i$-th row and $j$-th column of the matrix are denoted by $\vV_{i,*}$ and $\vV_{*,j}$, respectively.
For vectors $\vx$ and $\bbeta$ $\in$ $\mathbb{R}^{p}$, we denote their inner product $\langle \vx, \bbeta \rangle$ by $\vx(\bbeta)$, and the $k$-th entry of $\vx$ by $\vx(k)$. For two positive numbers $a,b$, we use $a\vee b$ and $a\wedge b$ to denote $\max\{a,b\}$ and  $\min\{a,b\}$ respectively.
Throughout the paper, we use $C$, $C'$, $C_1$ and $C_2$ to denote generic absolute constants, though the actual value may vary from case to case.
For two sequences $\{a_{n}\}$ and $\{b_{n}\}$, we denote $a_{n}\succeq b_{n}$ (resp. $a_{n}\preceq b_{n}$ ) if there exist positive constant $C$ (resp. $C'$) such that $a_{n} \geq Cb_{n}$ (resp. $a_{n} \leq C'b_{n}$). We denote $a_{n}\asymp b_{n}$ if both $a_{n}\succeq b_{n}$ and $a_{n}\preceq b_{n}$ hold. We denote $a_{n}\prec b_{n}$ (resp. $a_{n}\succ b_{n}$ ) if  $a_{n}=o(b_{n})$ (resp. $b_{n} =o(a_{n})$).
The $(1,\infty)$ norm and $(\infty, \infty)$ norm  of matrix A are defined by
$
\|A\|_{1,\infty}=\max_{1\leq j \leq p} \sum_{i=1}^{p}|A_{i,j}|$  and $\max_{1\leq i, j \leq n}\|A_{i,j}\|$ respectively.  For a finite subset $S$, we denote by $|S|$ its cardinality. We also write $A_{S,T}$ for the $|S|\times |T|$ sub-matrix with elements $(A_{i,j})_{i\in S, j\in T}$ and $A_{S}$ for $A_{S,S}$. For any squared matrix $A$, we define $\lambda_{min}(A)$ and $\lambda_{max}(A)$ as the smallest and largest eigenvalues of $A$, respectively.

\subsection{A brief review of the sliced inverse regression (SIR)}
SIR was first proposed by \citep{li1991sliced} to estimate the central space spanned by $\bbeta_{1},\ldots, \bbeta_{d}$ based on $n$ i.i.d. observations $(y_{i},\vx_{i})$, $i=1,\cdots,n$, from the multiple index model
$y=f(\bbeta_{1}^{\tau}\vx,...,\bbeta_{d}^{\tau} \vx,\epsilon)$, 
under the assumption  that $\vx$ follows an elliptical distribution and $\epsilon$ is Gaussian. SIR starts by
dividing the data into $H$ equal-sized slices according to the order statistics $y_{(i)}$. 
To ease notations and arguments, we assume that $n=cH$ and $\bbE[\vx]=0$, and  re-express the data as $y_{h,j}$ and $\vx_{h,j}$, where  $h$ refers to the slice number and $j$ refers to the order number within the slice, i.e.,
$y_{h,j}\leftarrow y_{(c(h-1)+j)}, \ \vx_{h,j} \leftarrow \vx_{(c(h-1)+j)}. $
Here $\vx_{(k)}$ is the concomitant of $y_{(k)}$.
Let the sample mean in the $h$-th slice 
be denoted by $ \overline{\vx}_{h,\cdot}$,
then  $\bLambda\triangleq var(\bbE[\vx|y])$ can be estimated by:
\begin{equation}\label{eqn:lambda}
\widehat{\bLambda}_{H}=\frac{1}{H}\sum_{h=1}^{H}\bar{\vx}_{h,\cdot}\bar{\vx}_{h,\cdot}^{\tau}=\frac{1}{H}\vX^{\tau}_{H}\vX_{H}
\end{equation}
where $\vX_{H}$  denotes the $p\times H$ matrix formed by the $H$ sample means, {\it i.e.}, $\vX_H=(\overline{\vx}_{1,\cdot},\ldots, \overline{\vx}_{H,\cdot})$.
Thus, $col(\bLambda)$ is estimated by ${col(\widehat{\vV}_{H})}$, where ${\widehat{\vV}_{H}}$ is the matrix formed by the $d$ eigenvectors associated to the largest $d$ eigenvalues of ${\widehat{\bLambda}_{H}}$.   
The ${col(\widehat{\vV}_{H})}$ is a consistent estimator of $col(\bLambda)$ under certain technical conditions \citep{duan1991slicing, hsing1992asymptotic, zhu2006sliced,li1991sliced, lin2017consistency}. 
It is shown in \cite{lin2017consistency,lin2018optimal} that, for single index models ($d=1$), $H$ can be chosen as a fixed number not depending on $\lambda(\bLambda)$, $n$, and $p$ for the asymptotic results to hold. Throughout this paper, we assume the following mild conditions. 
\begin{itemize}
\item[${\bf A1)}$]  $\vx \sim N(0,\bSigma)$ and there exist two positive constants $C_{\min}<C_{\max}$, such that $C_{\min}<\lambda_{\min}(\bSigma) \leq \lambda_{\max}(\bSigma)< C_{\max}$.
\end{itemize}
\begin{itemize}
\item
[${\bf A2)}$]   Matrix $var(\bbE[\vx|y])$ is non-vanishing, i.e., $\lambda_{\max}(var(\bbE[\vx|y])) >0$.
\end{itemize}

\subsection{Generalized Signal-to-Noise Ratio of Single Index Models 
}\label{subsec:problem}
We consider the following single index model:
\begin{equation}\label{eqn:sim}
y=f(\bbeta^\tau \vx, \bepsilon), \ \vx \sim N(0, \bSigma), \ \bepsilon\sim N(0, \sigma^2\vI),
\end{equation}
where $f(\cdot)$ is an unknown function. What we want to know is whether the coefficient vector $\bbeta$, when viewed as a whole, is zero. This can be formulated as a global testing problem as 
\[
H_0:\bbeta=0 \ \ \mbox{versus} \quad H_a: \bbeta\neq 0.
\]
When assuming the linear model $\vy =\bbeta^\tau\vx+\bepsilon$, whether we can separate the null and alternative depends on the interplay between $\sigma^2$ and the norm of $\bbeta$. More precisely, it depends on the signal-to-noise ratio (SNR) defined as
\begin{align*}
SNR=\frac{E[(\bbeta^{\tau}\vx)^{2}]}{\bbE[y^{2}]}=\frac{\|\bbeta\|_{2}^{2}\bbeta_{0}^{\tau}\bSigma\bbeta_{0}}{\sigma^{2}+\|\bbeta\|_{2}^{2}\bbeta_{0}^{\tau}\bSigma\bbeta_{0}}
\end{align*}
where $\bbeta_{0}=\bbeta/\|\bbeta\|_{2}$ \citep{janson2017eigenprism}.  Here $||\bbeta||_2$ is useful for benchmarking prediction accuracy for various model selection techniques such as AIC, BIC, or the Lasso. However, since there is an unknown link function $f(\cdot)$ in the single index model, the norm $||\bbeta||_2$ becomes non-identifiable. Without loss of generality, we restrict $||\bbeta||_2=1$ and have to find another quantity to describe the separability.

For the single index model (\ref{eqn:sim}),
to simplify the notation, 
let us use $\lambda$ to denote $\lambda_{\max}(var(\mathbb{E}[\vx|y]))$.
For linear models, we can easily show that
\begin{align*}
var(\bbE[\vx|y])=\frac{\bSigma\bbeta\bbeta^{\tau}\bSigma}{\bbeta_{0}^{\tau}\bSigma\bbeta_{0}\|\bbeta\|_{2}^{2}+\sigma^2} \mbox{ and } \lambda=\frac{\bbeta_{0}^{\tau}\bSigma\bSigma\bbeta_{0}\|\bbeta\|_{2}^{2}}{\bbeta_{0}^{\tau}\bSigma\bbeta_{0}\|\bbeta\|_{2}^{2}+\sigma^2}.
\end{align*}
Consequently, 
$
\lambda/SNR=\frac{\bbeta_{0}^{\tau}\bSigma\bSigma\bbeta_{0}}{\bbeta_{0}^{\tau}\bSigma\bbeta_{0}}
$.
When assuming condition A2), such a ratio is bounded by two finite limits. Thus, $\lambda$ can be treated as an equivalent quantity to the SNR for linear models, and is therefore named  as the generalized signal-to-noise ratio (gSNR) for single index models.

\begin{rmk}\normalfont  To the best of our knowledge, although SIR uses the estimation of $\lambda$ to determine the structural dimension (\cite{li1991sliced}), few investigations have been made towards theoretical properties of this procedure in high dimensions. The only work that uses $\lambda$ as a parameter to quantify the estimation error when estimating the direction of $\bbeta$ is \cite{lin2018optimal}, which, however, does not indicate explicitly what role $\lambda$ plays. 
 The aforementioned observation about $\lambda$ for single index models provides a useful intuition: $\lambda$ is a generalized notion of the SNR, and condition ${\bf A2)}$  merely requires that gSNR is non-zero.

\end{rmk}

\subsection{Global testing for single index models}

As we have discussed, \cite{arias2011global} and \cite{ingster2010detection} considered the 
testing problem \eqref{eqn:lr:detect}, which can be viewed as the determination of the detection boundary of gSNR. 
Through the whole paper, we consider the following testing problem:
\begin{equation}\label{eqn:raw:problem}
\begin{aligned}
\begin{cases}
H_{0}: &\bbeta_{1}=...=\bbeta_{p}=0,\\
H_{a}: \quad & \lambda (= gSNR) \quad \textrm{ is nonzero},
\end{cases}
\end{aligned}
\end{equation}
based on i.i.d. samples $\{(y_{i},\vx_{i}), i=1,\ldots, n\}$. Two models are considered: (i) the general single index model (SIM) defined in (\ref{eqn:sim}); and (ii) the single index model with additive noise (SIMa) defined as 
\begin{equation}\label{eqn:sima}
y=f(\bbeta^\tau \vx)+ \bepsilon, \ \vx \sim N(0, \bSigma), \ \bepsilon\sim N(0, \sigma^2\vI).
\end{equation}
We assume that conditions ${\bf A1)}$ and ${\bf A2)}$ hold for both models.

\section{The Optimal Test for Single Index Models} \label{sec:test}

\subsection{The detection boundary of linear regression}\label{subsec:lower:bound}
To set the goal and scope, we briefly review some related results on the detection boundary for linear models \citep{arias2011global,ingster2010detection}.

\begin{prop}\label{prop:low:bound:full:linear} Assume that $\vx_{i} \sim N(0,\vI_{p})$, $i=1,\cdots,n$, and that $\bbeta$ has at most $s$ non-zero entries. There is a test with both type I and II errors converging to zero for the testing problem in \eqref{eqn:lr:detect}  if and only if
\begin{align}
r^{2} \succ\frac{s\log(p)}{n}\wedge\frac{p^{1/2}}{n}\wedge\frac{1}{\sqrt{n}}.
\end{align}
\end{prop}

Assuming  $\vx\sim N(0,\vI_{p})$ and the variance of the noise is known, \cite{ingster2010detection} obtained the sharp detection boundary
(i.e., with exact asymptotic constant) for the above problem.
Since linear models are special cases of SIMa, which is a special subset of SIM, the following statement about the lower bound of detectability is a direct corollary of Proposition 
\ref{prop:low:bound:full:linear}.

\begin{cor}\label{cor:lower:bound1} 
$i)$ If $ s^{2}\log^2(p)\wedge p \prec n$,  then any test fails to separate the null and the alternative hypothesis asymptotically for SIM when
\begin{align}
\lambda \prec\frac{s\log(p)}{n}\wedge\frac{p^{1/2}}{n}.
\end{align} 
$ii)$ Any test fails to separate the null and the alternative hypothesis asymptotically for SIMa when
\begin{align}
\lambda \prec  \frac{s\log(p)}{n}\wedge\frac{p^{1/2}}{n} \wedge \frac{1}{\sqrt{n}}.
\end{align}
\end{cor}



\subsection{Single Index Models}
Moving from linear models to single index models is a big step. 
A natural and reasonable start is to consider tests based on the marginal correlation used for linear models \citep{ingster2010detection, arias2011global}. However, the following example shows that the marginal correlation fails for the single index models, indicating that we need to look for some other statistics to approximate the gSNR. 
\begin{example}\label{example:cubic}
Suppose that $\vx \sim N(0,\vI_{p})$, $\epsilon\sim N(0,1)$, and we have $n$ samples from the following model:
\begin{align}\label{eqn:example:model}
y=(x_{1}+...+x_{k})-(x_{1}+...+x_{k})^{3}/3k+\epsilon.
\end{align}
Simple calculation shows that $\bbE[\vx y]=0$. Thus, correlation-based methods do not work for this simple model.
On the other hand, since the link function $f(t)=t-t^3/3k$ is monotone when $|t|$ is sufficiently large, we know that $\bbE[\vx\mid y]$ is not a constant and $var(\bbE[\vx\mid y])\neq 0$.
\end{example}

Let $\lambda_0$ and $\lambda_0^a$ be two sequences such that
\[
\lambda_0\succ \frac{s\log(p)}{n}\wedge\frac{p^{1/2}}{n}, \ \ \lambda_0^a\succ  \frac{s\log(p)}{n}\wedge\frac{p^{1/2}}{n} \wedge \frac{1}{\sqrt{n}}.
\]
For a $p\times p$ symmetric matrix $A$ and a positive constant $k$ such that $ks<p$, we define 
\begin{align}\label{test:all:s}
\lambda^{(ks)}_{max}(A)=\max_{|S|=ks}\lambda_{max}(A_{S}).
\end{align}
For model $y=f(\bbeta^{\tau}\vx,\epsilon)$, in addition to the condition that $\lambda_0\prec \lambda$,  we further assume that  $s^{2}\log^2(p)\wedge p\prec n$.

 Let $\widehat{\bLambda}_{H}$ be the estimate of $var(\bbE[\vx|y])$ based on SIR. Let $\tau_n, \tau_n'$, and $\tau_n''$ be three quantities satisfying
 \begin{equation}\label{eqn:def:taus}
 \frac{\sqrt{p}}{n}\prec\tau_{n}\prec \lambda_0, \ \  \frac{s\log(p)}{n}\prec\tau_{n}'\prec \lambda_0, \ \ \frac{1}{\sqrt{n}}\prec \tau_{n}''\prec \lambda_0^a.
 \end{equation}
 
 We introduce the following two assistance tests
 \begin{itemize}
 \item[1.] Define 
 \[
 \psi_{1}(\tau_{n})=\bold{1}(\lambda_{\max}(\widehat{\bLambda}_H) >\frac{tr(\bSigma)}{n}+\tau_{n}).
 \]
 \item[2.] Define 
 \[ 
 \psi_{2}(\tau_{n}')=\bold{1}(\lambda^{(ks)}_{\max}(\widehat{\bLambda}_{H})>\tau_{n}'). 
 \]
 \end{itemize}
Finally, the {\bf S}pectral test {\bf S}tatistic based on {\bf S}IR, abbreviated as {\bf SSS}, is defined as
\begin{equation}\label{eqn:SSS}
\Psi_{SSS}=\max\{\psi_{1}(\tau_{n}),\psi_{2}(\tau_{n}')\}.
\end{equation}



To show the theoretical properties of SSS, we impose the following condition on the covariance matrix $\bSigma$.
\vspace*{1mm}
\begin{itemize}
\item[{\bf A3)}] There are at most $k$ non-zero entries in each row of $\bSigma$.
\end{itemize}

This assumption is first explicitly proposed in \cite{lin2017consistency}, which is partially motivated by the Separable After Screening (SAS) properties in \cite{ji2012ups}. 
In this paper, we assume such a relative strong condition and focus on establishing the detection boundary. 
This condition can be possibly relaxed by considering a larger class of covariance matrices 
\begin{align*}
\mathcal{S}(\gamma,\Delta)=\Big\{ |\bSigma_{jk}|\leq 1-(\log(p))^{-1}, \quad 
|\{k\mid  \bSigma_{jk}>\gamma\}|\leq \Delta
\Big\},
\end{align*}
which is used in \cite{arias2011detection} for analyzing linear models. 
Our condition contains $\mathcal{S}(0,\Delta)$ for some positive constant $\Delta$ and we could relax our constraint to some $\mathcal{S}(\gamma,\Delta)$. However, the technical details will be much more involved, which masks the importance of the main results. We thus leave it for a future investigation.

\begin{thm}\label{thm:main1}\normalfont
Assume that $s^{2}\log^2(p)\wedge p\prec n$, $\lambda\succ \lambda_0$,  and conditions ${\bf A1-A3)}$ hold.  
Two sequences $\tau_{n}$ and $\tau_{n}'$ satisfy the conditions in (\ref{eqn:def:taus}). Then, type I and type II errors of the test $\Psi_{SSS}(\tau_{n},\tau_{n}')$ converge to zero for the testing problem under SIM, i.e., we have
\begin{align*}
\bbE_{H_{0}}(\Psi_{SSS})+\bbE_{H_a}(1-\Psi_{SSS})\rightarrow 0.
\end{align*}
\end{thm}

Comparing with the test proposed in \cite{ingster2010detection}, our test statistics is a spectral statistics and depends on the first eigenvalue of $\widehat{\bLambda}_{H}$. It is adaptive in the moderate-sparsity scenario.
 In the high-sparsity scenario when $s^{2}\log^2(p)\prec p$, the SSS relies on $\psi_2(\tau_n')$, which depends on  the sparsity $s$ of the vector $\bbeta$. Therefore, SSS is not adaptive to the sparsity level.
Both \cite{arias2011detection} and \cite{ingster2010detection} introduced an (adaptive) asymptotically powerful test based on the higher criticism (HC) for the testing problem under linear models. It is an interesting research problem to develop an adaptive test using the idea of higher criticism for (\ref{eqn:raw:problem}).

\subsection{Optimal Test for SIMa}
When the noise  is assumed additive as in SIMa (\ref{eqn:sima}),
the detection boundary can be further improved. 
In addition to conditions {\bf A1-A3)}, $f$ is further assumed to satisfy the following condition: 
\begin{itemize}
\item[{\bf B)}] $f(z)$ is sub-Gaussian, $\bbE[f(z)]=0$ and $var(f(z) ) > C var(\bbE[z\mid f(z)+\epsilon])$ for some constant $C$, where $z, \epsilon \stackrel{iid}{\sim} N(0,1)$. 
\end{itemize}
Note that for any fixed function $f$ such that $var(\bbE[z\mid f(z)+\epsilon])\neq 0$, there exists a positive constant $C$ such that 
\begin{align}\label{temp:temp}
\frac{var(f(z))}{var(\bbE[z\mid f(z)+\epsilon])}>C.
\end{align}
By continuity, we know that \eqref{temp:temp} holds in a small neighbourhood of $f$, i.e., if $C$ is sufficiently small, condition $\vB)$ holds for a large class of functions.

First, we adopt the test $\Psi_{SSS}(\tau_{n},\tau_{n}')$ described in the previous subsection. 
Since the noise is additive, we include the ANOVA test,
\[
\psi_{3}(\tau_{n}'')=\bold{1}(t>\tau_{n}'')
\]
where $t=\frac{1}{n}\sum_{j=1}^{n}(y_{j}^{2}-1)$ and $\tau_{n}''$ is a sequence satisfying the condition (\ref{eqn:def:taus}). 
Combing this test with the test $\Psi_{SSS}(\tau_{n},\tau_{n}')$, we can introduce SSS assisted by ANOVA test (SSSa) as
\begin{align}\label{eqn:tilde_psi}
\Psi_{SSSa}(\tau_{n},\tau_{n}',\tau_{n}'')=\max\{\Psi_{SSS}(\tau_{n},\tau_{n}'), \psi_{3}(\tau_{n}'')\}
\end{align}

We then have the following result.

\begin{thm}\label{thm:main2}\normalfont
Assume that $\lambda \succ \lambda_0^a$ and the conditions ${\bf A1-A3)}$ and ${\bf B)}$ hold.  
Assume that the sequences $\tau_{n}$, $\tau_{n}'$ and $\tau_n''$ satisfy condition (\ref{eqn:def:taus}), then  type I and type II errors of the test ${\Psi}_{SSSa}(\tau_{n},\tau_{n}',\tau_{n}'')$ converge to zero for the testing problem under SIMa, i.e., we have
\begin{align*}
\bbE_{H_{0}}({\Psi}_{SSSa})+\bbE_{H_{a}}(1-{\Psi}_{SSSa})\rightarrow 0.
\end{align*}
\end{thm}

{\it 
{\bf Example continued.} 
For the example in (\ref{eqn:example:model}),
we calculated the test statistic $\psi_{SSS}$ defined by  (\ref{eqn:SSS}) under both the null and alternative hypotheses. Figure \ref{fig:example:sir} shows the histograms of such a statistic under both hypotheses, demonstrating  a perfect separation between the null and alternative.  For this example, $\lambda_{\max}^{ks}(\widehat{\bLambda}_H)$ has more discrimination power than $\lambda_{\max}(\widehat{\bLambda}_H)$.

\begin{figure}
\centering
\includegraphics[height=75mm, width=75mm]{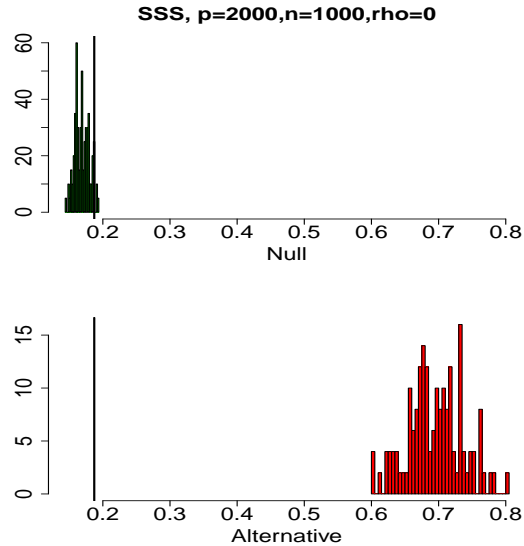}
\caption{ The histograms of $\lambda_{max}^{ks}(\hat{\bLambda}_H)$ for the model (\ref{eqn:example:model}). The top panel corresponds to the scores under the null and the bottom one corresponds to the scores under the alternative. The "black" vertical line is the 95\% quantile under the null.
}\label{fig:example:sir}
\end{figure}

}

\subsection{Computationally efficient test}
Although the test $\Psi_{SSS}$ (and ${\Psi}_{SSSa}$) is rate optimal, it is computationally inefficient. Here we propose an efficient algorithm to approximate $\lambda_{max}^{(ks)}(\widehat{\bLambda}_H)$ via a convex relaxation, which is similar to the convex relaxation method for estimating the top eigenvector of a semi-definite  matrix  in \cite{adamczak2008smallest}. To be precise, given the SIR estimate $\widehat{\bLambda}_{H}$ of $var(\bbE[\vx \mid y])$, consider the following semi-definite programming (SDP) problem:
\begin{equation}
\begin{aligned}
\widetilde{\lambda}_{\max}^{(ks)}(\widehat{\bLambda}_H) \triangleq \max \quad &tr (\widehat{\bLambda}_{H}\vM)\\
\mbox{subject to} \quad &tr(\vM)=1, \quad |\vM|_{1}\leq ks,\\
 &\vM  \mbox{ is semi-definite positive}.
\end{aligned}
\end{equation} 

With $\widetilde{\lambda}^{(ks)}_{\max}(\widehat{\bLambda}_H)$, for a sequence $\tau_{n}'$ satisfying the condition in (\ref{eqn:def:taus}), i.e., $\frac{s\log(p)}{n}\prec\tau_{n}'\prec \lambda_0$, a computationally feasible test is
 \[ \widetilde{\psi}_{2}(\tau_{n}')=\bold{1}(\widetilde{\lambda}^{(ks)}_{\max}(\widehat{\bLambda}_{H})>\tau_{n}').\]
 
Then, for any sequence $\tau_{n}$ satisfying the inequality in (\ref{eqn:def:taus}), we define the following computationally feasible alternative of $\Psi_{SSS}$:
	\begin{align}
\widetilde{\Psi}_{SSS}=\max\{\psi_{1}(\tau_{n}),{\widetilde\psi}_{2}(\tau_{n}')\}
\end{align}
\begin{thm}\label{thm:main3}\normalfont
Assume that $s^{2}\log^2(p)\wedge p\prec n$, $\lambda\succ\lambda_0$ and conditions ${\bf A1-A3)}$ hold. Then, type I and type II errors of the test $\widetilde{\Psi}_{SSS}(\tau_{n},\tau_{n}')$ converge to zero for the testing problem under SIMa, i.e., we have
\begin{align*}
\bbE_{H_{0}}(\widetilde{\Psi}_{SSS})+\bbE_{H_{a}}(1-\widetilde{\Psi}_{SSS})\rightarrow 0.
\end{align*}
\end{thm}

Similarly, if we introduce the test
\begin{align}
\widetilde{\Psi}_{SSSa}(\tau_{n},\tau_{n}',\tau_{n}'')=\max\{\widetilde{\Psi}_{SSS},\psi_{3}(\tau_{n}'')\},
\end{align}
for three sequences $\tau_{n}$,$\tau_{n}'$ and $\tau_{n}''$, then we have
\begin{thm}\label{thm:main4}\normalfont
Assume that $ \lambda\succ \lambda_0^a$ and  conditions ${\bf A1-A3)}$ and ${\bf B)}$ hold.  
The test $\widetilde{\Psi}_{SSSa}(\tau_{n},\tau_{n}',\tau_{n}'')$ is asymptotically powerful for the testing problem under SIMa, i.e., we have
\begin{align*}
\bbE_{H_{0}}({\widetilde{\Psi}}_{SSSa})+\bbE_{H_{a}}(1-{\widetilde{\Psi}}_{SSSa})\rightarrow 0.
\end{align*}
\end{thm}


Theorem \ref{thm:main2} and Theorem \ref{thm:main4} not only establish the detection boundary of gSNR for single index models, but also open a door of thorough understanding of semi-parametric regression with a Gaussian design. 
It is shown in \cite{lin2018optimal} that for single index models satisfying conditions {\bf A1)}, {\bf A2)}, {\bf A3)}, one has
\begin{align}
\sup_{\widehat{\beta}} \inf_{\mathfrak{m}} \bbE_{\mathfrak{m}}\|P_{\widehat{\beta}}-P_{\beta}\|^{2}_{F} \asymp 1\wedge \frac{s\log(ep/s)}{n\lambda}.
\end{align}
This implies that
the necessary and sufficient condition for obtaining a consistent estimate of the projection operator $P_{\bbeta}$ is $\frac{s\log(ep/s)}{n}\prec \lambda$. 
On the other hand, Theorems \ref{thm:main2} and \ref{thm:main4} state that, for single index models with additive noise, if $\frac{s\log(p)}{n}\wedge\frac{p^{1/2}}{n}\wedge\frac{1}{\sqrt{n}}\prec\lambda$, then one can detect the existence of gSNR (a.k.a. non-trivial direction $\bbeta$). 
Our results thus imply for SIMa that,
if $\frac{p^{1/2}}{n}\wedge \frac{1}{\sqrt{n}}\prec\lambda\prec\frac{s\log(p)}{n}$, one can detect the existence of non-zero $\bbeta$, but can not provide a consistent estimation of its direction. To estimate the location of non-zero coefficient, we must tolerate a certain error rate such as the false discovery rate (\cite{Benjamini:Hochberg:1995}). For example, the knockoff procedure (\cite{barber2015controlling}), SLOPE (\cite{Su:Candes:2016}), and UPT(\cite{Ji:Zhao:2014}) might be extended to  single index models.


\subsection{Practical Issues}\label{sec:practical}

In practice, we do not know whether the noise is additive or not. Therefore, we only consider the test statistic $\widetilde{\Psi}_{SSS}$. 
Condition (\ref{eqn:def:taus}) provides us a theoretical basis for choosing the sequences $\tau_n$ and $\tau_n'$. In practice, however, we determine these thresholds by simulating the null distribution of $\lambda_{\max}(\widehat{\bLambda}_H)$ and $\tilde{\lambda}_{max}^{(ks)}(\widehat{\bLambda}_H)$. Our final algorithm is as follows.

\begin{algorithm}[H]
\caption{Spectral test Statistic based on SIR (SSS) Algorithm\label{alg:SSS}}
\begin{algorithmic}
\vspace*{1mm}
\\ 1.  Calculate $\lambda_{\max}(\widehat{\bLambda}_H)$ and $\widetilde{\lambda}_{\max}^{(ks)}(\widehat{\bLambda}_H)$ for the 
given input $(\vx, \vy)$;
\vspace*{1mm}
\\ 2. Generate $\vz=(z_1,\cdots,z_n)$, where $z_i\iid N(0,1)$;
\vspace*{1mm}
\\ 3. Calculate $\lambda_{\max}(\widehat{\bLambda}_H)$ and $\widetilde{\lambda}_{\max}^{(ks)}(\widehat{\bLambda}_H)$ based on $(\vx, \vz)$;
\vspace*{1mm}
\\ 4. Repeat Steps 2 and 3 $N (=100)$ times to get two sequences of $\lambda_{max}$ and $\widetilde{\lambda}_{\max}^{(ks)}$. Let $\tau_n$ and $\tau_n'$ be the 95\% quantile of these two simulated sequences;
\vspace*{1mm}
\\ 5. Reject the null if $\lambda_{\max}(\widehat{\bLambda}_H) >\tau_n$ and/or $\tilde{\lambda}_{max}^{(ks)}(\widehat{\bLambda}_H)>\tau_n'$. 
\end{algorithmic}
\end{algorithm}

\section{Numerical Studies}\label{sec:studies}

 Let $\bbeta$ be the vector of coefficients and let $\mathcal{S}$ be the active set, $\mathcal{S}=\{i:\beta_i\neq 0\}$, for which we simulated $\beta_i \stackrel{iid}{\sim} N(0,1)$.
Let $\vx$ be the random design matrix with each row following $N(0, \bSigma)$. 
We consider two types of covariance matrices: (i) $\bSigma = (\sigma_{ij})$ with $\sigma_{ii}=1$ and $\sigma_{ij}=\rho^{|i-j|}$; 
and (ii) $\sigma_{ii}=1$, $\sigma_{ij}=\rho$ when $i,j\in\mathcal{S}$ or $i,j\in \mathcal{S}^c$, and $\sigma_{ij}=\sigma_{ji}=0.1$ when $i\in\mathcal{S}, j\in\mathcal{S}^c$. 
The first one represents a covariance matrix which is essentially sparse and we choose $\rho$ among 0, 0.3, 0.5, and 0.8. 
The second one represents a dense covariance matrix with $\rho$ chosen as 0.2. 
In all the simulations, $n=1,000$, $p$ varies among 100, 500, 1,000, and 2,000 and the number of replication is 100.
The random error $\bepsilon$ follows $N(0,\vI_{n})$.
We consider the following models:

\begin{itemize}
\item[I.] $\vy= 0.02*\left( 16 \vx\bbeta - exp(\vx\bbeta) \right) + \bepsilon$, where $|\mathcal{S}|=7$;
\item[II.]  $\vy=  0.2 * \sin(\vx\bbeta/2)*exp(\vx\bbeta/2) + \epsilon $, where $|\mathcal{S}|=10$;
\item[III.] $\vy= 0.8*\left(\vx\bbeta -  (\vx\bbeta)^3/15\right) + \epsilon $, where $|\mathcal{S}|=5$;
\item[IV.] $\vy = \sin(\vx\bbeta)*exp(\vx\bbeta/10)*\epsilon$, where $|\mathcal{S}|=10$.
\end{itemize}

\begin{table}[!htbp]
\caption{Power comparison of SSS and HC for four models I-IV for different parameter settings. Symbol ``$\ast$" indicates the type (ii) covariance matrix.}\label{tab:power:1}
\begin{tabular}{|c|c|c|c|c|c|c|c|c|c|}
\hline
Model & Dim &$\rho$& SSS & HC & Model & Dim &$\rho$& SSS & HC \\
\hline
\multirow{20}{*}{I}  & \multirow{5}{*}{100} & 0 & 1.00 & 0.16 & \multirow{20}{*}{II} & \multirow{5}{*}{100} & 0 & 0.98 & 0.12 \\
& & 0.3 & 1.00 & 0.29 & & & 0.3 & 0.97 & 0.16 \\ 
& & 0.5 & 0.99 & 0.54 & & & 0.5 & 0.96 & 0.24 \\
& & 0.8 & 1.00 & 0.93 & & & 0.8 & 1.00 & 0.37\\
& & 0.2$^\ast$ & 0.90 & 0.35 & & & 0.2$^\ast$ & 0.96 & 0.56 \\
\cline{2-5}\cline{7-10}
& \multirow{5}{*}{500} & 0 & 0.98 & 0.16 & & \multirow{5}{*}{500} & 0 & 0.87 & 0.06\\
& & 0.3 & 0.99 & 0.18 & & & 0.3 & 0.80 & 0.09 \\
& & 0.5 & 0.97 & 0.34 & & & 0.5 & 0.82 & 0.13 \\
& & 0.8 & 0.98 & 0.71 & & & 0.8 & 0.83 & 0.14 \\
& & 0.2$^\ast$ & 0.52 & 0.25 & & & 0.2$^\ast$ & 0.77 & 0.32\\
\cline{2-5}\cline{7-10}
& \multirow{5}{*}{1000} & 0 & 0.89 & 0.19 & & \multirow{5}{*}{1000} & 0 & 0.81 & 0.09\\
& & 0.3 & 0.88 & 0.16 & & & 0.3 & 0.74 & 0.06 \\
& & 0.5 & 0.91 & 0.33 & & & 0.5 & 0.77 & 0.08 \\
& & 0.8 & 0.96 & 0.53 & & & 0.8 & 0.84 & 0.11 \\
& & 0.2$^\ast$ & 0.37 & 0.30 & & & 0.2$^\ast$ & 0.69 & 0.25\\
\cline{2-5}\cline{7-10}
& \multirow{5}{*}{2000}&0 &0.92 & 0.18 & & \multirow{5}{*}{2000} & 0 & 0.75 & 0.11\\
& & 0.3 & 0.86 & 0.25 & & & 0.3 & 0.68 & 0.12 \\
& & 0.5 & 0.83 & 0.43 & & & 0.5 & 0.68 & 0.13 \\
& & 0.8 & 0.90 & 0.60 & & & 0.8 & 0.81 & 0.10 \\
& & 0.2$^\ast$ & 0.43 & 0.17 & & & 0.2$^\ast$ & 0.63 & 0.41 \\
\hline
\hline
\hline
\multirow{20}{*}{III}  & \multirow{5}{*}{100} & 0 & 1.00 & 0.21 & \multirow{20}{*}{IV} & \multirow{5}{*}{100} & 0 & 0.89 & 0.01 \\
& & 0.3 & 1.00 & 0.25 & & & 0.3 & 0.91 & 0.03 \\ 
& & 0.5 & 1.00 & 0.63 & & & 0.5 & 0.89 & 0.04 \\
& & 0.8 & 1.00 & 1.00 & & & 0.8 & 1.00 & 0.10\\
& & 0.2$^\ast$ & 0.98 & 0.78 & & & 0.2$^\ast$ & 0.94 & 0.07 \\
\cline{2-5}\cline{7-10}
& \multirow{5}{*}{500} & 0 & 0.99 & 0.11 & & \multirow{5}{*}{500} & 0 & 0.70 & 0.03\\
& & 0.3 & 1.00 & 0.12 & & & 0.3 & 0.57 & 0.04 \\
& & 0.5 & 0.98 & 0.11 & & & 0.5 & 0.57 & 0.07 \\
& & 0.8 & 0.99 & 0.22 & & & 0.8 & 0.69 & 0.09 \\
& & 0.2$^\ast$ & 0.62 & 0.72 & & & 0.2$^\ast$ & 0.45 & 0.08\\
\cline{2-5}\cline{7-10}
& \multirow{5}{*}{1000} & 0 & 0.99 & 0.11 & & \multirow{5}{*}{1000} & 0 & 0.55 & 0.07\\
& & 0.3 & 0.97 & 0.06 & & & 0.3 & 0.56 & 0.04 \\
& & 0.5 & 0.97 & 0.18 & & & 0.5 & 0.51 & 0.09 \\
& & 0.8 & 0.92 & 0.10 & & & 0.8 & 0.73 & 0.06 \\
& & 0.2$^\ast$ & 0.60 & 0.59 & & & 0.2$^\ast$& 0.44 & 0.08\\
\cline{2-5}\cline{7-10}
& \multirow{5}{*}{2000}&0 & 0.96 & 0.16 & & \multirow{5}{*}{2000} & 0 & 0.58 & 0.07\\
& & 0.3 & 0.97 & 0.19 &  & & 0.3 & 0.47 & 0.07 \\
& & 0.5 & 0.93 & 0.15 & & & 0.5 & 0.45 & 0.09 \\
& & 0.8 & 0.88 & 0.10 & & & 0.8 & 0.61 & 0.02 \\
& & 0.2$^\ast$ & 0.59 & 0.58 & & & 0.2$^\ast$ & 0.40 & 0.08 \\
\hline

\end{tabular}
\end{table}

We choose $k=1$ in the Algorithm \ref{alg:SSS}.
If we calculate $N(=100)$ test statistics for each replication, it will take an extremely long time. Therefore, in the simulation, we calculate $\tau_n$ and $\tau_n'$ slightly different from  Algorithm \ref{alg:SSS}. For each generated data set, we simulated only one vector $\vz$ where $\vz\sim N(0,\vI_n)$ and calculate the statistic $\lambda_{\max}(\widehat{\bLambda}_H)$ and $\widetilde{\lambda}_{\max}^{(ks)}(\widehat{\bLambda}_H)$. The $\tau_n$ and $\tau_n'$ are chosen as 95\% quantile from the corresponding sequence for all the replications.

For each generated data, we also calculated Cor-HC scores according to \cite{arias2012detection}. The threshold $c_{hc}$ is chosen according to the same scheme as choosing the thresholds $\tau_n$ and $\tau_n'$. Namely, we calculated the Cor-HC scores based on $\vz$ where $\vz\sim N(0,\vI_n)$. The threshold $c_{hc}$ is the 95\% quantile of these simulated scores.
The hypotheses is rejected if the Cor-HC score is greater than $c_{hc}$. The power for both methods is calculated as the average number of rejections out of 100 replications. These numbers are reported in Tables \ref{tab:power:1}.

\begin{figure}
\includegraphics[width=50mm, height=50mm]{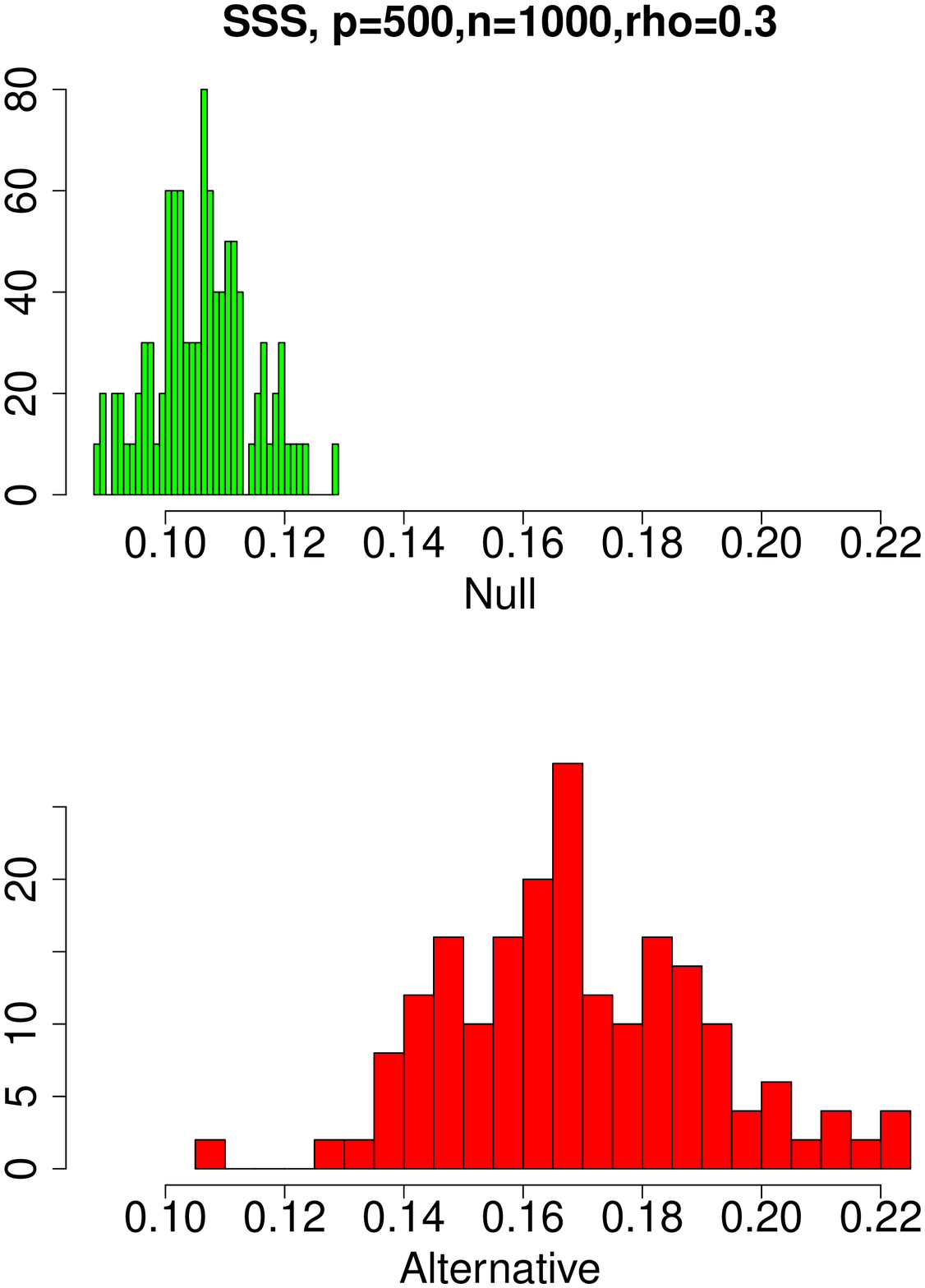}
\includegraphics[width=50mm, height=50mm]{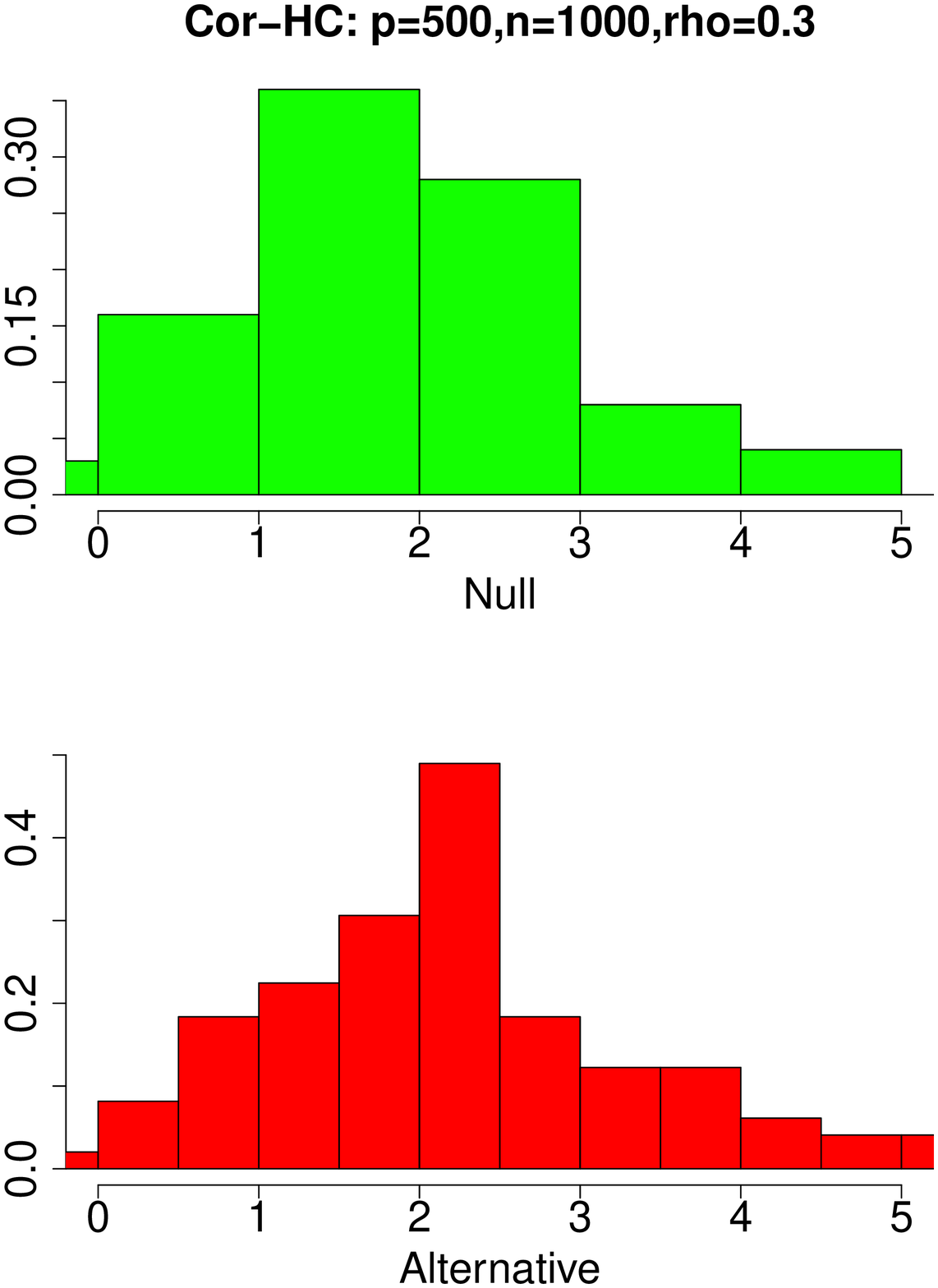}
\caption{Model III, $n=1,000$, $p=500$, type (i) covariance matrix, $\rho=0.3$.}\label{fig:setting4}
\end{figure}

It is clearly seen that the power of SSS decreases when the dimension $p$ increases. Nevertheless, the power of SSS is better than the one based on Cor-HC except for one case. 
In Figure \ref{fig:setting4}, we plot the histogram of the statistic $\widetilde{\lambda}_{\max}^{(ks)}(\widehat{\bLambda}_H)$ under the null in the top-left panel and the histogram of this statistic under the alternative in the bottom-left panel for Model III when $p=500$ and $\rho=0.3$ for type (i) covariance matrix. It is clearly seen that the test statistic $\Psi_{SSS}$ are well separated under the null and alternative. However, Cor-HC fails to distinguish between the null and alternative as shown in the two panels on the right side.  

To see how the performance of Cor-HC varies, we consider the following model 
\begin{itemize}
\item[V.] $\vy= \kappa \vx\bbeta - exp(\vx\bbeta) + \bepsilon $,  where $|\mathcal{S}|=7,\kappa=1,3,5,\cdots, 19.$
\end{itemize}

\begin{figure}[!h]
\centering
 \includegraphics[height=75mm, width=75mm]{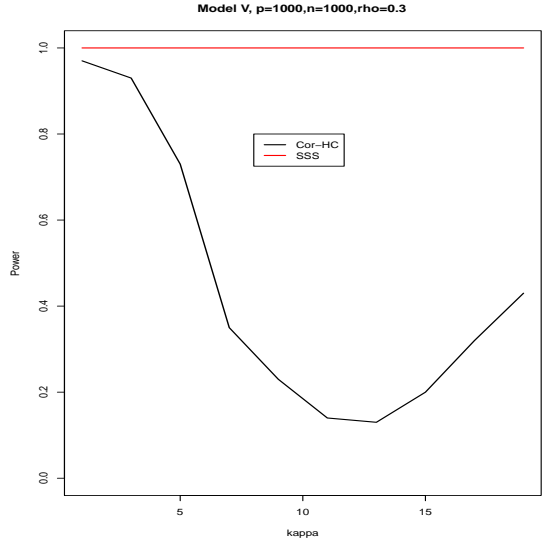}
\caption{Power: Model V, $n=1,000, p=1,000, \rho=0.3$ for type (i) covariance matrix.}\label{fig:vary:coef}
\end{figure}

Set $n=1,000$, $p=1,000$, and $\rho=0.3$ for  type (i) covariance matrix, and the power of both methods are displayed in Figure \ref{fig:vary:coef}. The coefficient $\kappa$ determines the magnitude of the marginal correlation between the active predictors and the response. It is seen that when $\kappa$ is close to 16, representing the case of diminishing marginal correlation, the power of Cor-HC dropped to the lowest. Under all the models, SSS is more powerful in detecting the existence of the signal.

To observe the influence of the signal-to-noise ratio on the power of the tests, we consider the following two models
\begin{itemize}
\item[VI.] $\vy= ( 15 \vx\bbeta - exp(\vx\bbeta) )*\kappa + 4\bepsilon $, where $|\mathcal{S}|=7$;
\item[VII.] $\vy = \sin(\vx\bbeta)*exp(10\vx\bbeta\kappa)*\epsilon$, where $|\mathcal{S}|=10$.
\end{itemize}
Here $\kappa=0.01,0.02,\ldots, 0.10$.

\begin{figure}[!h]
\centering
\includegraphics[height=60mm, width=60mm]{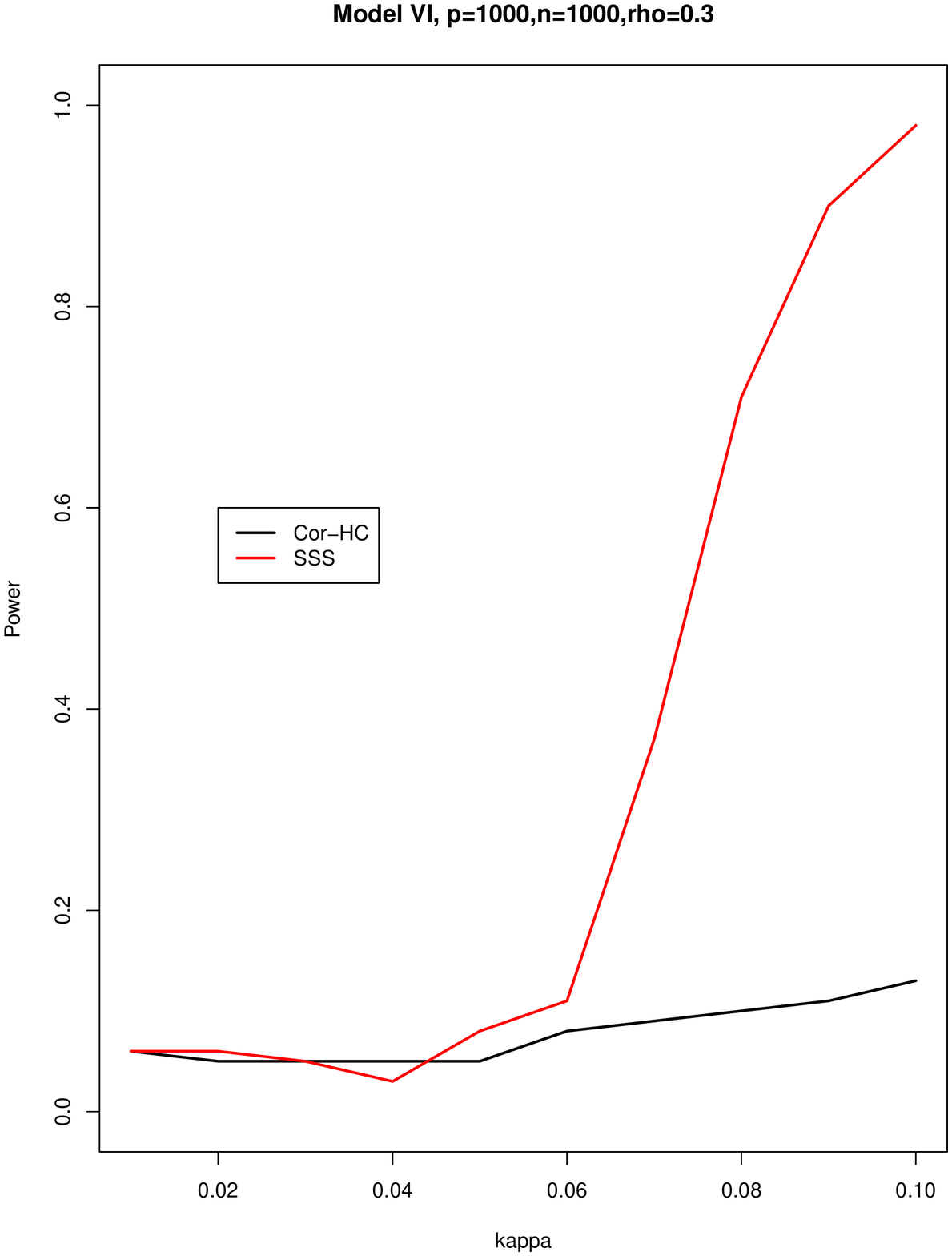}
\includegraphics[height=60mm, width=60mm]{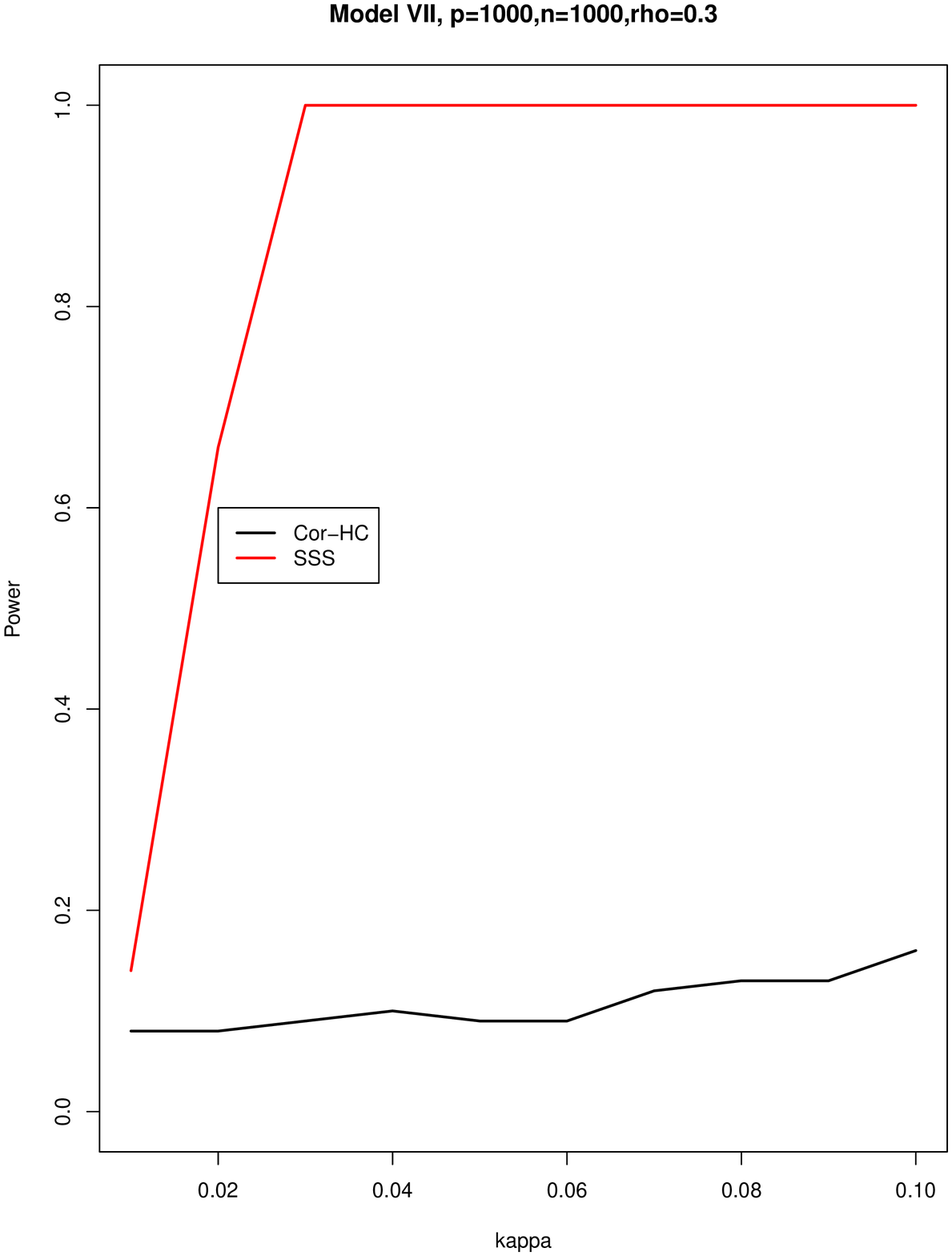}
\caption{Power: Models VI and VII, $n=1,000, p=1,000, \rho=0.3$ for the type (i) covariance matrix.}\label{fig:vary:SNR}
\end{figure}

Set $n=1,000, p=1,000$ and $\rho=0.3$, we plot the power of both methods against the coefficient $\kappa$ in Figure \ref{fig:vary:SNR}. It is clearly seen that for both examples there is a sharp "phase-transition" for the power of SSS as the signal strength increases, validating our theory about the detection boundary. In both examples SSS is much more powerful than Cor-HC.







\section{Discussion}\label{sec:discuss}
Assuming that $var(\bbE[\vx\mid y])$ is non-vanishing, we show in this paper that $\lambda$, the unique non-zero eigenvalue of $var(\bbE[\vx\mid y])$ associated with the single index model, is a generalization of the SNR. We demonstrate a surprising similarity between linear regression and single index models with Gaussian design: the detection boundary of gSNR for the testing problem (\ref{eqn:raw:problem}) under SIMa matches that of SNR for linear models (\ref{eqn:lr:detect}).
This similarity provides an additional support to the speculation that ``the rich theories developed for linear regression can be extended to the single/multiple index models" \citep{lin2018sparse,chen1998can}.

Besides  the gap we explicitly depicted between detection and estimation boundaries, we provide here several other directions which might be of interests to researchers.
First, although  this paper only deals with single index models, the results obtained here
are very likely extendable to multiple index models. 
Assume that the noise is additive and let $0<\lambda_{d}\leq ...\leq \lambda_{1}$ be the non-zero eigenvalues associated with the matrix $var(\bbE[\vx|y])$ of a multiple index model.  
Similar arguments can show that the $i$-th direction is detectable if $\lambda_{i} \succ \frac{\sqrt{p}}{n}\wedge \frac{s\log(p)}{n}\wedge\frac{1}{\sqrt{n}}$.  New  thoughts 
and technical preparations might be needed for a rigorous argument for determining the lower bound of the detection boundary.  
Second, the framework can be extended to study  theoretical properties of  other sufficient dimension reduction algorithms such as SAVE and directional regression (\cite{lin2017consistency,lin2018optimal,lin2018sparse}).


\section{Acknowledgment}
We  thank Dr. Zhisu Zhu for his generous help with SDP.

\vspace{0.3in}

\begin{center}
{\bf APPENDIX: PROOFS}
\end{center}

\begin{appendices}

\section{Assisting Lemmas}
Since our approaches are based on the technical tools developed in \cite{lin2017consistency,lin2018optimal,lin2018sparse},  we briefly recollect the necessary (modified) statements without proofs below.

\begin{lem}
Let $z_{j}\iid N(0,1), j=1,\ldots,p$. Let $\sigma_1,\ldots,\sigma_p$ be $p$ positive constants satisfying $\sigma_{1}\leq \ldots\leq \sigma_{p}$. Then for any $0<\alpha\leq \frac{1}{\sigma_{p}^{2}}\sum_{j}\sigma_{j}^{4}$, we have
\begin{align}
\bbP\left(\sum_{j} \sigma^{2}_{j}\left( z_{j}^{2}-1 \right) >\alpha \right)\leq \exp\left(-\frac{\alpha^{2}}{4\sum\sigma^{4}_{j}}\right).
\end{align}
\end{lem}

\begin{lem}\label{lem:2}
Suppose that a $p\times H$ matrix $\vX$ formed by $H$ $i.i.d.$ $p$ dimensional vector $\vx_{j}\sim N(0,\bSigma)$, $j=1,\ldots,H$ where $0<C_{1}\leq \lambda_{\min}(\bSigma)\leq \lambda_{\max}(\bSigma) \leq C_{2}$ for some constants $C_{1}$ and $C_{2}$. We have
\begin{align}
\|\frac{1}{p}\vX^{\tau}\vX-\frac{tr(\bSigma)}{p}\vI_{H}\|_{F}>\alpha
\end{align}
with probability at most $4H^{2}\exp\left(-\frac{Cp\alpha^{2}}{H^{2}}\right)$ for some positive constant $C$.
In particular, we know that
\begin{align}
\lambda_{\max}\left( \vX\vX^{\tau} /p\right) =\lambda_{\max}\left( \vX^{\tau}\vX/p \right) \leq tr(\bSigma)/p+\alpha
\end{align} 
happens with probability at least 
$1-4H^{2}\exp\left(-\frac{Cp\alpha^{2}}{H^{2}}\right)$.

\end{lem}

\begin{lem}\label{lem:3}  Assume that $p^{1/2}\prec n\lambda$. 
Let $\vM=\left(\begin{array}{cc}
B_{1} & 0\\
B_{2} & B_{3}\\
0 & B_{4} 
\end{array} \right)$ be a $p\times H$ matrix, where $B_{1}$ and $B_{2}$ are scalar, $B_{3}$ is a $1\times (H-1)$ vector and $B_{4}$ is a $(p-2)\times(H-1)$ matrix satisfying 
\begin{equation}
\begin{aligned}
& (1-\frac{1}{2\nu})\lambda \leq B_{1}^{2}\leq (1+\frac{1}{2\nu})\lambda  \\
& \left\|\left(\begin{array}{cc}
B_{2}^{2} & B_{2}B_{3} \\
B_{3}^{\tau}B_{2} & B_{3}^{\tau}B_{3}+B_{4}^{\tau}B_{4}
\end{array} \right)-\frac{A}{n}\vI_{H}\right\|_{F}\leq \frac{\sqrt{p}\alpha}{n}
\end{aligned}.
\end{equation}
for a constant $\nu>1$  where $\alpha\prec \frac{n\lambda}{p^{1/2}}$. Then we have
\begin{align}
\lambda_{\max}\left(\vM\vM^{\tau}\right) >\frac{A}{n}-\frac{\sqrt{p}\alpha}{n}+(1-\frac{1}{2\nu})\lambda.
\end{align}
\end{lem}

\paragraph{Sliced approximation inequality}\label{append:sliced stable} The next result is referred to as `key lemma' in \cite{lin2017consistency,lin2018optimal,lin2018sparse}
, which depends on the following sliced stable condition. 
\begin{definition}[Sliced stable condition]\label{def:sliced:stable} 
For  $0<\mathfrak{a}_{1} < 1<\mathfrak{a}_{2} $, 
 let $\mathcal{A}_H(\mathfrak{a}_{1},\mathfrak{a}_{2})$ denote all partitions $\{ -\infty = a_0 \leq a_2 \leq \ldots \leq a_{H} = +\infty\}$ of $\mathbb{R}$ satisfying that
 \[
 \frac{\mathfrak{a}_{1}}{H} \leq \mathbb{P}(a_h \leq Y \leq a_{h + 1}) \leq \frac{\mathfrak{a}_{2}}{H}.
 \]
A curve $m(y)$ is $\vartheta$-sliced stable with respect to y, if there exist positive constants $\mathfrak{a}_{1},\mathfrak{a}_{2},\mathfrak{a}_{3}$ and large enough $H_{0}$ such that for any $H>H_{0}$, for any partition in $\mathcal{A}_H(\mathfrak{a}_{1},\mathfrak{a}_{2})$ and any $\gamma \in \mathbb{R}^{p}$ , one has:
\begin{align}\label{sliced_inequality}
\frac{1}{H}\sum_{h=1}^{H}var\left(\gamma^{\tau}m(y)\big| a_{h-1}\leq y < a_{h}\right)\leq \frac{\mathfrak{a}_{3}}{H^{\vartheta}}  var\left(\gamma^{\tau}m(y)\right).
\end{align}
A curve is sliced stable if it is $\vartheta$-sliced stable for some positive constant $\vartheta$.
\end{definition}

The sliced stable condition is a mild condition. \cite{neykov2015support} derived the sliced stable condition from a modification of the regularity condition proposed in \cite{hsing1992asymptotic}. The inequality \eqref{sliced_inequality} implies the following deviation inequality for multiple index models. For our purpose, we modified it for  single index models.

\begin{lem}\label{lem:4}
Assume that Conditions ${\bf A1)}$, ${\bf A2)}$ and the sliced stable condition (for some $\vartheta>0$) hold  in the single index model $y=f(\bbeta^{\tau}\vx, \epsilon)$. Let $\widehat{\bLambda}_{H}$ be the SIR estimate of $\bLambda=var(\bbE[\vx \mid y])$, and let $P_{\bLambda}$ be the projection matrix associated with the column space of $\bLambda$.
For any vector $\beta\in \mathbb{R}^{p}$ and any $\nu>1$, 
let $E_{\beta}(\nu)=\Big\{~ \Big| \beta^{\tau}\left(P_{\bLambda}\widehat{\bLambda}_{H}P_{\bLambda}-\bLambda\right)\beta\Big| \leq \frac{1}{2\nu}\beta^{\tau}\bLambda\beta\Big\}$.
There exist positive constants $C_{1},C_{2}$,$C_{3}$ and $C_{4}$ such that
for any $\nu>1$ and $H$ satisfying that $H^{\vartheta}>C_{4}\nu$, one has  
\begin{align}
\mathbb{P}\left(  \bigcap_{\beta}E_{\beta} \right)\geq 1-C_{1}\exp\left(-C_{2}\frac{n\lambda_{\max}(\bLambda)}{H\nu^{2}} +C_{3}\log(H)\right).
\end{align} 
\end{lem}
\cite{lin2018stability} recently proved a similar deviation inequality without the sliced stable condition. 

\section{Proof of Theorems}
 \vspace{4mm}
\paragraph{Proof of Theorem \ref{thm:main1}} Theorem \ref{thm:main1} follows from the following Lemma \ref{lem:5} and Lemma \ref{lem:6}.
\begin{lem}\label{lem:5}
 Assume that $p^{1/2}\prec n\lambda_0$,
and $\tau_{n}$ be a sequence such that $\frac{\sqrt{p}}{n}\prec \tau_{n} \prec \lambda_0$. Then, as $n\rightarrow \infty$, we have: 

\vspace*{2mm}
\hspace{3mm}  $i)$ Under $H_{0}$, i.e., if $y \independent \vx$, then $\lambda_{\max}(\widehat{\bLambda}_H)<\frac{tr(\bSigma)}{n}+\tau_{n}$ with probability converging to 1;

\vspace*{2mm}
\hspace{3mm} $ii)$ Under $H_{1}$, if $\lambda\succ \lambda_0$, then $\lambda_{\max}(\widehat{\bLambda}_H)>\frac{tr(\bSigma)}{n}+\tau_{n}$ with probability converging to 1. 
\proof 
$i)$ If $y \independent \vx$, we know that $\frac{1}{\sqrt{H}}\bSigma^{-1/2}\vX_{H}$ is a $p\times H $ matrix with entries $i.i.d.$ to $N(0,\frac{1}{n})$. From Lemma \ref{lem:2}, we know that
\begin{align}
\lambda_{\max}(\frac{1}{H}\vX_{H}^{\tau}\vX_{H})\leq \frac{tr(\bSigma)}{n}+\tau_{n}
\end{align}
with probability at least $1-4H^{2}\exp\left(-\frac{Cn^{2}\tau_{n}^{2}}{H^{2}p} \right)$ 
which $\rightarrow 1$ as $n \rightarrow \infty$. 

$ii)$ For any event $\omega$,  there exist $p\times p$ orthogonal matrix $S$ and $H\times H$ orthogonal matrix $T$  such that 
\begin{align}
S\vX_{H}(\omega)T=\left(\begin{array}{cc}
Z_{1} & 0\\
Z_{2} & Z_{3}\\
0 & Z_{4}
\end{array} \right)
\end{align}
where $Z_{1}$, $Z_{2}$ are two scalars, $Z_{3}$ is a $1\times(H-1)$ vector and $Z_{4}$ is a $(p-2)\times (H-1)$ matrix. Lemma \ref{lem:4} and Lemma \ref{lem:2} imply  that there exist a constant $A$ and an events set $\Omega$, such that $\bbP\left(\Omega^{c}\right) \rightarrow 0$ as $n\rightarrow \infty$.
For any $\omega \in \Omega$, one has
\begin{equation}\label{Condition:C}
\begin{aligned}
& (1-\frac{1}{2\nu})\lambda \leq Z_{1}^{2}\leq (1+\frac{1}{2\nu})\lambda  \\
& \left\|\left(\begin{array}{cc}
Z_{2}^{\tau}Z_{2} & Z_{2}^{\tau}Z_{3} \\
Z_{3}^{\tau}Z_{2} & Z_{3}^{\tau}Z_{3}+Z_{4}^{\tau}Z_{4}
\end{array} \right)-\frac{tr(\bSigma)}{n}\vI_{H}\right\|_{F}\leq \frac{\sqrt{p}\alpha}{n}.
\end{aligned}
\end{equation}
Lemma \ref{lem:3} implies that 
\begin{align}
\lambda_{\max}(\frac{1}{H}\vX_{H}\vX_{H})^{\tau} \geq \frac{tr(\bSigma)}{n} -\frac{\sqrt{p}}{n}\alpha+(1-\frac{1}{2\nu})\lambda \succ \frac{tr(\bSigma)}{n}+\tau_{n}.
\end{align}
\end{lem}

\begin{lem}\label{lem:6} Assume that $\frac{s\log(p)}{n}\prec\lambda_0$. 
Let $\tau_{n}$ be a sequence such that $\frac{s\log(p)}{n}\prec \tau_{n} \prec \lambda_0$. Then, as $n\rightarrow \infty$, we have: 

\hspace{3mm} i) if $y\independent \vx$, then $\lambda^{(ks)}_{\max}\left(\widehat{\bLambda}_{H}\right)<\tau_{n}$ with probability converging to 1;

\hspace{3mm} ii) if $\lambda \succ \lambda_0$, then 
$\lambda^{(ks)}_{\max}\left(\widehat{\bLambda}_{H}\right)>\tau_{n}$ with probability converging to 1. 
\proof 
 $i)$ If $y \independent \vx$, we know that $\frac{1}{\sqrt{H}}\vE_{H}=\frac{1}{\sqrt{H}}\bSigma^{-1/2}\vX_{H}$ is a $p\times H $ matrix with entries $i.i.d.$ to $N(0,\frac{1}{n})$. Thus
\begin{align*}
\lambda^{(ks)}_{\max}\left(\widehat{\bLambda}_{H}\right)=\lambda^{(ks)}_{max}\left(\frac{1}{H}\bSigma^{1/2}\vE_{H}\vE_{H}^{\tau}\bSigma^{1/2} \right) \mbox{ and }  \lambda_{\max}^{(ks)}\left(\vD^{1/2}\vE_{H}\vE_{H}^{\tau}\vD^{1/2}\right).
\end{align*}
are identically distributed
where $\vD$ is diagonal matrix consisting of the eigenvalues of $\bSigma$. For any subset $S\subset [p]$, let $\vX_{H,S}=\vD_{S}^{1/2}\vE_{S,H}$ where $\vE_{S,H}$ is a submatrix of $\vE_{H}$ consisting of the rows in $S$. Note that 
\begin{align}
\lambda_{\max}\left(\left(\vD^{1/2}\vE_{H}\vE_{H}^{\tau}\vD^{1/2} \right)_{S}\right)
=\lambda_{\max}\left(\frac{1}{H}\vX_{H,S}\vX_{H,S}^{\tau}\right).
\end{align}
Thus, by Lemma \ref{lem:3}, we have
\begin{align}
\lambda_{\max}\left(\frac{1}{H}\vX_{H,S}\vX_{H,S}^{\tau}\right)<tr(\vD_{S})/n+ \alpha  \leq \frac{ks\lambda_{\max}(\bSigma)}{n}+ \alpha
\end{align}
with probability at least $1-4H^{2}\exp\left(-\frac{Cn^{2}\alpha^{2}}{H^{2}s} \right)$. Let $\alpha=C\frac{s\log(p)}{n}$ for some sufficiently large constant $C$. 
Since ${p \choose ks} \leq \left(\frac{ep}{ks}\right)^{ks}$, we know that $\lambda^{(ks)}_{\max}(\widehat{\bLambda}_{H})\leq C\frac{s\log(p)}{n} \prec \tau_{n}$ with probability converges to 1.

$ii)$ Let $\boldeta$ be the eigenvector associated to the largest eigenvalue of $\bLambda$. Thus $|supp(\boldeta)|=ks$. From Lemma \ref{lem:4}, we know that 
\begin{align}
\widehat{\lambda}_{\max}^{(ks)}(\widehat{\bLambda}_{H})
\geq \boldeta^{\tau}\widehat{\bLambda}_{H}\boldeta
\geq (1-\frac{1}{2\nu})\lambda 
\end{align}
with probability converges to 1. 
Thus, $\lambda^{(ks)}_{\max}\left(\widehat{\bLambda}_{H} \right)>\tau_{n}$ with probability converges to 1. 
\end{lem}
\epf

 \vspace{4mm}
 \paragraph{Proof of Theorem \ref{thm:main2}} Theorem \ref{thm:main2} follows from the Theorem \ref{thm:main1} and the following Lemma \ref{lem:7}.

\begin{lem}\label{lem:7}Assume that $\frac{1}{\sqrt{n}}\prec\lambda_0^a$. Let $\tau_{n}$ be a sequence such that $\frac{1}{\sqrt{n}}\prec\tau_{n}\prec \lambda_0^a$, $\tau_{n}\rightarrow 0$. Then we have

\hspace*{3mm}  $i)$ If $y\independent \vx$, then $t<\tau_{n}$ with probability converging to 1.

\hspace*{3mm} $ii)$ If $\lambda\succ\lambda_0^a$, then $t>\tau_{n}$ with probability converging to 1.
\proof
(i) Since $y\independent \vx$ ,we know that $\bbE[t]=0$. Let $z_{j}=y_{j}^{2}-1$, then we have
\begin{align}
\bbP\left(\frac{1}{n}\sum_{j}z_{j}>\tau_{n} \right)\leq \exp\left(-Cn\tau_{n}^{2} \right)
\end{align}
for some constant $C$. In other words, the probability of $t>\tau_{n}$  converges to 0 as $n \rightarrow \infty$.

(ii) If $\lambda\succ\lambda_0^a$, we have $var(f(z)) \geq C\lambda$ and $\bbE[y^{2}-1]\geq C\lambda$ for some constant $C$ . 
Let $z_{j}=y_{j}^{2}-1$, $j=1,\ldots,n$. Since $f(x_{j})$, $j=1,\ldots,n$ are sub-Gaussian, we know that
\begin{align}
\bbP\left( \frac{1}{n}\sum_{j}z_{j}>\bbE[y^{2}-1]+\delta \right)\leq \exp\left(-Cn\delta^{2} \right)
\end{align}
By choosing $\delta=C\bbE[y^{2}-1]$ for some constant $C$, we know that the probability of   
$
t\geq  (C+1)\lambda\succ \tau_{n}    
$
converges to 1. 
\end{lem}
\epf

 \vspace{4mm}
\paragraph{Proof of Theorem \ref{thm:main3}}  Theorem \ref{thm:main3} and \ref{thm:main4} follows from the following Lemma, the Theorem \ref{thm:main1} and the Theorem \ref{thm:main2}.

\begin{lem}\label{lem:8}
Assume that $\frac{s\log(p)}{n}\prec \lambda_0$. 
Let $\tau_{n}$ be a sequence such that $\frac{s\log(p)}{n} \prec \tau_{n}\prec \lambda_0$. Then we have:

\hspace{3mm} $i)$ if $y\independent \vx$, then $\widetilde{\lambda}^{(ks)}_{\max}\left(\widehat{\bLambda}_{H}\right)<\tau_{n}$ with probability converging to 1;

\hspace{3mm} $ii)$ if $\lambda\succ\lambda_0$, then $\widetilde{\lambda}^{(ks)}_{\max}\left(\widehat{\bLambda}_{H}\right)>\tau_{n}$ with probability converging to 1. 
\proof $i)$ 
Under $H_{0}$, i.e., $y\independent \vx$, the entries of $\frac{1}{\sqrt{H}}\bSigma^{-1/2}\vX_{H}$ is identically distributed as $N(0,\frac{1}{n})$. Thus, if $1 \prec \alpha \prec \frac{n\tau_{n}}{s\log(p)}$, we have 
\begin{align}
\max_{(i,j)} \left| \widehat{\bLambda}_{H}(i,j) \right| \leq \frac{\alpha\log(p)}{n}
\end{align}
with probability at least $1-p^{2}\exp\left(-C\alpha^{2}\log(p)^{2} \right)$ for some constant $C$ which converges to 1 as $n \rightarrow \infty$. Since (See e.g., Lemma 6.1 in \cite{berthet2013optimal})
\begin{align}
\widetilde{\lambda}^{(ks)}_{\max}\left( \widehat{\bLambda}_{H}\right) \leq \lambda_{\max}\left(st_{\frac{\alpha\log(p)}{n}}\left(\widehat{\bLambda}_{H} \right)\right)+ks\frac{\alpha\log(p)}{n}\prec\tau_{n}
\end{align} 
where $st_{z}(A)_{i,j}=sign(A_{i,j})(A_{i,j}-z)_{+}$, 
we know that $i)$ holds.

$ii)$ Follows from that $\widetilde{\lambda}^{(ks)}_{max}(\widehat{\bLambda}_{H}) \geq \lambda^{(ks)}_{\max}(\widehat{\bLambda}_{H})$. \epf
\end{lem}

\end{appendices}

\bibliographystyle{plainnat}
\bibliography{sir}

\end{document}